\documentclass{amsart}

\usepackage{epsfig}
\usepackage{amsmath}
\usepackage{amssymb}
\usepackage{amscd}
\usepackage{graphicx}
\usepackage{xcolor}
\usepackage{float}
\usepackage{verbatim}
\usepackage{bbm}
\usepackage{cancel}
\usepackage{enumerate}

\newtheorem{thm}{Theorem}[section]
\newtheorem{lem}[thm]{Lemma}
\newtheorem{prop}[thm]{Proposition}

\newtheorem{ques}[thm]{Question}

\newtheorem{defn}[thm]{Definition}

\theoremstyle{remark}
\newtheorem{rem}[thm]{Remark}

\def \N {\mathbb N}

\def \B {\mathcal B}

\def \Z {\mathbb Z}

\def \Q {\mathcal Q}

\def \P {\mathcal P}
\def \S {\mathcal S}
\def \M {\mathcal M}

\def \xt {$(X,T)$}

\def \eps {\varepsilon}

\def \sq {sequence}

\def \ys {$(Y,S)$}
\def \tl {topological}
\def \im {invariant measure}
\def \inv {invariant}

\def \ds {dynamical system}

\def \diam {\mathsf{diam}}

\numberwithin{equation}{section}

\begin{document}

\title{Lifting generic points}

\author{Tomasz Downarowicz and Benjamin Weiss}

\address{\vskip 2pt \hskip -12pt Tomasz Downarowicz}

\address{\hskip -12pt Faculty of Pure and Applied Mathematics, Wroc\l aw University of Technology, Wroc\l aw, Poland}

\email{downar@pwr.edu.pl}

\medskip
\address{\vskip 2pt \hskip -12pt Benjamin Weiss}

\address{\hskip -12pt Einstein Institute of Mathematics,
The Hebrew University of Jerusalem}

\email{weiss@math.huji.ac.il}

\begin{abstract}
Let \xt\ and \ys\ be two \tl\ \ds s, where \xt\ has the weak specification property.
Let $\xi$ be an \im\ on the product system $(X\times Y, T\times S)$ with marginals
$\mu$ on $X$ and $\nu$ on $Y$, with $\mu$ ergodic. Let $y\in Y$ be quasi-generic for $\nu$.
Then there exists a point $x\in X$ generic for $\mu$ such that the pair $(x,y)$ is quasi-generic for $\xi$. This is a generalization of a similar theorem by T.\ Kamae, in which \xt\ and \ys\ are full shifts on finite alphabets.
\end{abstract}

\maketitle

\section{Introduction}\label{S1}

\noindent \emph{\ \ All terminology used freely in the introduction is explained in the preliminaries.}
\medskip

Let  $\pi:(X,T) \rightarrow (Y,S)$ be an extension of a compact dynamical system $(Y,S)$ and suppose that $\nu$ is an ergodic measure for $S$. This measure can always be lifted to an invariant measure on $X$ (by the Hahn-Banach theorem). It then follows that there exists an ergodic measure $\mu$ that projects to $\nu$. Clearly, any generic point for $\mu$ will project to a generic point for $\nu$. It is natural to ask whether all $\nu$-generic points can be obtained in this way. In other words: does every $\nu$-generic point have a $\mu$-generic lift? It is not difficult to show that if $\mu$ is a unique lift of $\nu$ then the answer is yes. In fact, in this case if $y \in Y$ is generic for $\nu$ then any $x\in \pi^{-1}(y)$ will be generic for $\nu$. However, if the extension of $\nu$ is not unique then the answer might be negative. Such examples can be obtained as follows:

Consider a minimal almost 1-1 extension $\pi:X \rightarrow Y$ where $Y$ is strictly ergodic but $X$ is not (cf.\ Furstenberg--Weiss \cite{FW} for examples of such systems). Then all \im s on $X$ project to the unique measure $\nu$ on $Y$. In this situation all points in $Y$ are generic for $\nu$. Now, let $y\in Y$ be a point with a unique preimage $x\in X$ (by assumption such a point exists). Clearly, $x$ can be generic for at most one measure on $X$, hence there will be a measure $\mu$ on $X$ (which is a lift of $\nu$) such that $y$ (generic for $\nu$) does not lift to a point generic for $\mu$. 

However, this example does not provide an answer to a more subtle question: Does every $\nu$-generic point have a generic lift (without specifying for which measure extending $\nu$)? In general, the answer to such a relaxed question is also negative. We will show this using the same example as before and the following theorem:

\begin{thm}\label{one} Let $(X,T)$ be a \tl\ \ds\ with (at least) two different ergodic measures $\mu$ and $\nu$, both having full \tl\ support. Then there is a  dense $G_{\delta}$-set $B\subset X$ and a continuous function $f$ on $X$ such that for any $x\in B$ the ergodic averages 
$$
A_n(f,x)=\frac1n\sum_{i=0}^{n-1}f(T^ix) 
$$
oscillate.
\end{thm}

\begin{proof} Since the measures differ, there exists a continuous function $f$ on $X$ whose
integral with respect to $\mu$ is greater than one while its integral with respect to $\nu$ is less than zero. Now, for a natural number $N$ we define
$$
E_N = \Bigl\{ x \in X :\exists_{n > N}\  \frac{1}n\sum_{i=0}^{n-1} f(T^ix) > 1\Bigr\}.
$$
This set is clearly open and it is dense since the generic points for $\mu$ are dense. Define a similar set $F_N$ replacing ``$>1$'' by ``$<0$''. Then the the desired set $B$ is  
the countable intersection 
$$
B=\bigcap_{N\ge1}\bigl(E_N\cap F_N).
$$
By the Baire theorem, this set is a dense $G_\delta$, and clearly no $x\in B$ is generic for any measure.
\end{proof}

Now let us go back to the example. Since the system $(X,T)$ is minimal, all its \im s have full \tl\ support. By Theorem~\ref{one}, there is a dense $G_\delta$-set $B$ of points which are not generic for any measure. As in any minimal almost 1-1 extension, the ``singleton fibers'' (that is points which are unique preimages of their images) also form a dense $G_\delta$-set (call it $A$) in $X$. Then the intersection $A\cap B$ is nonempty and any point in its image is generic (for $\nu$) but has no generic~lift.

\medskip
Before we discuss a positive result we need to mention two important issues. The first one is the phenomenon of \emph{quasi-generating} \im s, that is, generating them along a subsequence of averages. Replacing the term ``generic'' by ``quasi-generic'' may lead to either stronger or weaker results, depending on where the replacement is done (in the assumption or in the thesis). The second issue is a specific way of extending a system by \emph{joining it} with another system. Practically, any extension can be viewed as a joining of the system with its extension (the joining is then supported by the graph of the factor map) but it is often essential to know that the extension can be obtained as a joining with a system having some specific properties (such as ergodicity, specification property, etc.) which the entire extension does not necessarily enjoy.
\smallskip

In the early 70's, Teturo Kamae studied normal \sq s and the phenomenon of normality preserving sub\sq s. In symbolic dynamics a \sq\ over a finite alphabet is \emph{normal} if it is generic for the uniform Bernoulli measure. An increasing sub\sq\ of natural numbers $y=(n_k)_{k\ge1}$ \emph{preserves normality} if $x|_y=(x_{n_k})_{k\ge1}$ is normal for any normal \sq\ $(x_n)_{n\ge1}$. Few years earlier, B.\ Weiss \cite{W} proved that sub\sq s of positive lower density which are completely deterministic preserve normality. A sub\sq\ $y$ is \emph{completely deterministic} if its indicator function, viewed as an element of the shift on two symbols, quasi-generates only measures of entropy zero. Kamae \cite{K} proved the opposite implication: \emph{only} completely deterministic sub\sq s preserve normality. Given a non-deterministic sub\sq\ $y$ (i.e., one that quasi-generates some measure $\nu$ of positive entropy), he needed to find a normal \sq\ $x$ such that $x|_y$ is not normal. Skipping the details, let us just say that he needed to ``pair'' the sub\sq\ $y$ with a normal (i.e., generic for the uniform Bernoulli measure~$\lambda$) \sq\ $x$, such that the pair $(x,y)$ is generic for a specific joining $\xi$ of $\lambda$ and~$\nu$. In order to do so, he proved a more general theorem, which motivates our current work. We take the liberty of rephrasing the statement in the language that we use throughout this note. 

\begin{thm}\label{kamae} {\rm\cite{K}} Let $\xi$ be a joining of two \im s, $\mu$ and $\nu$, supported on symbolic systems $\Lambda_1^\N$ and $\Lambda_2^\N$, respectively ($\Lambda_1$ and $\Lambda_2$ are finite alphabets). Let $y\in\Lambda_2^\N$ be quasi-generic for $\nu$, i.e., generates $\nu$ along a sub\sq\ of averages indexed by $\mathcal J=(n_k)_{k\ge1}$. Then there exists $x\in\Lambda_1^\N$ such that the pair $(x,y)$ generates $\xi$ along~$\mathcal J$. If $\mu$ is ergodic then $x$ can be chosen generic for $\mu$.
\end{thm} 

This theorem found another application in the work of G.\ Rauzy \cite{R}, who studied normality preservation in a different meaning. Let us identify all real numbers with their expansions in some fixed base $b\ge2$. A real number is called \emph{normal (in base $b$)} if its expansion is a normal \sq. A real number $y$ \emph{preserves normality} if $x+y$ is normal for any normal number~$x$. Rauzy proved that a number $y$ preserves normality if and only if the expansion of $y$ is completely deterministic.

\smallskip

Notice that Theorem \ref{kamae} is actually very strong. First of all, it applies to any situation when a ``symbolic'' measure $\nu$ is lifted to a ``symbolic'' measure $\xi$. Also note that $\nu$ is not assumed ergodic, it suffices that it admits a quasi-generic point~$y$ (which is always true within a full shift). If $\mathcal J=\N$ then $y$ is simply generic for~$\nu$ and the theorem allows to lift it to a pair $(x,y)$ generic for $\xi$. Even when $\mathcal J$ is an essential sub\sq\ (and there is no hope to make the lift $(x,y)$ generic), as soon as $\mu$ is ergodic, the point $x$ ``paired'' with $y$ still can be generic rather than just quasi-generic. The only weakness of the theorem is that $x$ is found within the \emph{full shift} $\Lambda_1^\N$, even when $\mu$ is supported by a proper subshift. In other words, the theorem does not allow to lift $y$ within an \emph{a priori} given \tl\ (symbolic) extension.

And our paper focuses exactly on this problem. Our goal is to find conditions under which the ``paired'' point $x$ (generic for $\mu$) can be found within the \emph{a priori} given \tl\ system $(X,T)$ being joined with $(Y,S)$. The conditions turn out to be: ergodicity of $\mu$ (like in the original theorem), and the \emph{weak specification property} of $(X,T)$. We prove the following:
\begin{thm}\label{lgp}
Let \xt\ and \ys\ be \tl\ \ds s and let $\xi$ be an \im\ on the product system $(X\times Y,T\times S)$ with marginals $\mu$ and $\nu$ on $X$ and $Y$, respectively. Assume that the system \xt\ has the weak specification property and that $\mu$ is ergodic under $T$. Suppose also that $y\in Y$ is quasi-generic for the measure $\nu$. Then there exists a point $x\in X$, generic for $\mu$, such that the pair $(x,y)$ is quasi-generic for $\xi$. 
\end{thm}

Let us mention that the weak specification property is satisfied by many systems such as ergodic mixing Markov shifts, ergodic toral automorphisms, and in fact any endomorphisms of compact Abelian groups for which the Haar measure is ergodic (see \cite{D}). An advantage of our result is that it is not restricted to symbolic systems and that $x$ is found within the space~$X$. A disadvantage is that the pair $(x,y)$ is only quasi-generic for $\xi$, even when $\mathcal J=\N$. Theorem \ref{lgp} will be applied in the forthcoming paper \cite{B-D}, where Rauzy's equivalence between normality preservation and determinism is generalized to a wider context. Namely, addressed is the following problem:

\begin{ques}\label{ques} 
Let $T:X\to X$ be a surjective endomorphism of a compact metrizable Abelian group, such that the Haar measure $\lambda$ on $X$ is $T$-ergodic and has finite entropy. Let us call a point $x\in X$ \emph{normal} if it is generic for $\lambda$.
Is it true that $y$ \emph{preserves normality} (i.e., if $x+y$ is normal for any normal $x\in X$) if and only if $y$ is completely deterministic?
\end{ques}

In \cite{B-D} we prove sufficiency relatively easily, but the harder direction (necessity) is shown only for selected groups $X$ (tori, solenoids, and countable direct products $\bigoplus_{n\ge1}\Z_d$ ($d\ge2$)). In all these cases Theorem \ref{lgp} plays a crucial role in the proofs. The necessity in full generality remains open.

\medskip
Our paper is organized as follows. Section \ref{S2} contains all necessary definitions and notational conventions. In Section \ref{S3} we provide three key lemmas together with auxiliary propositions needed in their proofs. The propositions are quite standard while the lemmas may be considered of independent interest. Finally, in Section \ref{S4} we present the proof of Theorem \ref{lgp}.

\section{Preliminaries}\label{S2}

Let \xt\ be a \tl\ \ds, where $X$ is a compact metric space and $T$ is a continuous surjection. By $\M(X)$ we will denote the set of all Borel probability measures on $X$. Since no other measures will be considered, the elements of $\M(X)$ will be shortly called \emph{measures}. By $\M_T(X)$ we will denote the subset of $\M(X)$ containing all measures that are $T$-\inv, i.e., such that $\mu(T^{-1}A)=\mu(A)$ for all Borel sets $A\subset X$. When the transformation $T$ is fixed, the elements of $\M_T(X)$ will be called \emph{\im s}. The sets $\M(X)$ and $\M_T(X)$ are equipped with the weak* topology, which makes both these sets compact convex and metrizable with a convex metric.\footnote{By definition, a \sq\ $(\mu_n)_{n\ge1}$ of measures converges in the weak* topology to a measure~$\mu$ if, for any continuous (real or complex) function $f$ on $X$, the integrals $\int f\,d\mu_n$ converge to $\int f\,d\mu$. One of standard convex metrics compatible with this topology is given by
$$
d(\mu,\nu)=\sum_{n=1}^\infty 2^{-n}\Bigl|\int f_n\,d\mu - \int f_n\,d\nu\Bigr|,
$$
where $(f_n)_{n\ge1}$ is a \sq\ of continuous functions on $X$ with values in the interval $[0,1]$, linearly dense in the space $C(X)$ of all continuous real functions on $X$.} 
 It is well known that the extreme points of $\M_T(X)$ are precisely the \emph{ergodic measures}, i.e., \im s $\mu$ such that $\mu(A\,\triangle\, T^{-1}A)=0 \implies \mu(A)\in\{0,1\}$, for any Borel set $A\subset X$.

We will be using the following notation. For two integers $a\le b$, by $[a,b]$ we will denote the interval of integers $\{a,a+1,a+2,\dots,b\}$. Given a point $x\in X$ and $0\le a\le b$, by $x[a,b]$ we denote the ordered finite segment of the orbit of $x$:
$$
x[a,b]=(T^ax,T^{a+1}x,\dots,T^bx),
$$
while by $\mu_{x[a,b]}$ we will understand the normalized counting measure on $x[a,b]$:
$$
\mu_{x[a,b]}=\frac1{b-a+1}\sum_{n=a}^b\delta_{T^nx}
$$ 
(here $\delta_x$ denotes the Dirac measure concentrated at $x$). We will call this measure the \emph{empirical measure} associated with the orbit segment. 


A point $x$ is said to \emph{quasi-generate} (or be \emph{quasi-generic} for) a measure $\mu$ if $\mu$ is an accumulation point of the \sq\ of measures $(\mu_{x[0,n]})_{n\ge1}$ (any such measure $\mu$ is \inv). In this case there exists an increasing \sq\ of natural numbers $\mathcal J=(n_k)_{k\ge1}$ such that $\lim_{k\to\infty}\mu_{x[0,n_k]}=\mu$. We will say that $x$ \emph{generates $\mu$ along~$\mathcal J$}. If the \sq\ $(\mu_{x[0,n]})_{n\ge1}$ converges to $\mu$ then we say that $x$ \emph{generates} (or \emph{is generic} for) $\mu$ (in other words, generic = generic along $\N$). It follows from the Pointwise Ergodic Theorem that every ergodic measure $\mu$ possesses generic points (in fact $\mu$-almost all points are such). The follwing obvious fact will be used several times:

\begin{rem}\label{eqseq}
Two increasing \sq s of natural numbers, say $(n_k)_{k\ge1}$ and $(m_k)_{k\ge1}$ will be called \emph{equivalent} if
$\lim_{k\to\infty}\frac{n_k}{m_k}=1$. It is obvious that the upper (and lower) densities of any subset of $\N$ evaluated along equivalent \sq s are the same. If a point $x$ generates a measure $\mu$ along a \sq\ $(n_k)_{k\ge1}$ then it generates $\mu$ along any \sq\ $(m_k)_{k\ge1}$ equivalent to $(n_k)_{k\ge1}$. 
\end{rem}

Other key notions in this paper are those of a specification and the specification property:

\begin{defn}\label{spec}
\mbox{}
\begin{enumerate}
\item Consider a (finite or infinite) \sq\ of nonnegative integers:
\begin{align*}
a_1\le &\,b_1<a_2\le b_2<\dots<a_{N_1}\le b_{N_1}, \ \text{where } N_1\in\N, \text{ or }\\
a_1\le &\,b_1<a_2\le b_2<a_3\le b_3<\dots\,.
\end{align*}
Let $D=\bigcup_N[a_N,b_N]$ (where $N$ ranges over either $[1,N_1]$ or $\N$). 
By a \emph{specification} with domain $D$ we will mean any function 
$$
\S:D\to X 
$$ 
such that for each $N$ there exists a point $x_N$ such that for each $n\in[a_N,b_N]$ we have
$$
\mathcal S(n) = T^n(x_N).
$$
Since $T$ is surjective, we can equivalently demand that $\S(n)=T^{n-a_N}(x_N)$. 
\item By $\S[a_N,b_N]$ we mean the ordered tuple $(\S(a_N),\S(a_N+1),\dots,\S(b_N))$ which equals $x_N[a_N,b_N]$ (or $x_N[0,b_N-a_N])$. 
\item The numbers $l_N=b_N-a_N+1$ and $g_N=a_{N+1}-b_N-1$ will be called the \emph{orbit segment lengths} and \emph{gaps} of the specification, respectively.
\item If $D$ is finite then the \emph{empirical measure} associated with $\mathcal S$ is defined as
$$
\mu_{\mathcal S}=\frac1{|D|}\sum_{n\in D}\delta_{\mathcal S(n)}.
$$
\item An infinite specification $\S$ \emph{generates} (or \emph{is generic} for) a measure $\mu$ along a \sq\ $\mathcal J=(n_k)_{k\ge1}$ if the measures associated to the specification $\S$ restricted to $D\cap[0,n_k]$ converge to $\mu$ (in general, $\mu$ need not be \inv).
\item We say that (the orbit of) a point $x\in X$ $\eps$-shadows the specification $\S$ if 
$$
\forall_{n\in D}\ d(\S(n),T^n(x))<\eps.
$$

\item A system \xt\ has the \emph{weak specification property} if for every $\eps>0$ there exists a function $M_\eps:\N\to\N$ satisfying $\lim_{l\to\infty}\frac{M_\eps(l)}l=0$, such that any finite specification (with any finite number $N_1$ of orbit segments) satisfying, for each $N\in[1,N_1]$, the inequality $g_N\ge M_\eps(l_{N+1})$ is $\eps$-shadowed by an orbit.
\end{enumerate}
\end{defn}

The last condition asserts, roughly speaking, that any appropriately spaced finite sequence of orbit segments (where each gap is adjusted to the length of the following segment, according to the function $M_\eps$) can be $\eps$-shadowed by a single orbit. 


\section{Preparatory statements}\label{S3}

The proof of Theorem \ref{lgp} relies on three key lemmas. The first one is concerned with increasingly good shadowing of certain infinite specifications.

\begin{lem}\label{goodshad}
Let \xt\ be a \tl\ \ds\ with the weak specification property with a family of fuctions $\{M_\eps:\eps>0\}$. Let $(\eps_k)_{k\ge0}$ be a summable \sq\ of positive numbers. Let $D=\sum_{N=1}^\infty[a_N,b_N]$ and let $\S:D\to X$ be an infinite specification 
satisfying, for some increasing \sq\ $(N_k)_{k\ge0}$ of nonnegative integers starting with $N_0=0$, the following condition: 
For each $k\ge1$ and all $N\in[N_{k-1}+1,N_k]$ we have
\begin{equation}\label{szu}
g_N\ge M_{\eps_k}(l_{N+1}). 
\end{equation}
Then there exists a point $x_0$ such that
\begin{equation}\label{szer0}
\lim_{n\in D}\ d(\S(n),T^nx_0)=0.
\end{equation}
\end{lem}

\begin{proof}
Let $\S_1$ denote the specification $\S$ restricted to the initial $N_1$ orbit segments. This finite specification satisfies the inequality $g_N\ge M_{\eps_1}(l_{N+1})$, hence it can be $\eps_1$-shadowed by the orbit of some point $x_1\in X$. 

We continue by induction. Suppose that we have found a point $x_k\in X$ which satisfies
\begin{equation}\label{szer}
\forall_{i\in[1,k]}\ \forall_{N\in[N_{i-1}+1,N_i]}\ \forall_{n\in[a_N,b_N]}\ d(\S(n),T^nx_k)\le\sum_{j=i}^k\eps_j.
\end{equation}
We define a new specification $\S_{k+1}$ on 
$$
[0,b_{N_k}]\cup\bigcup_{N=N_k+1}^{N_{k+1}}[a_N,b_N]
$$ 
as follows: 
We let $\S_{k+1}[0,b_{N_k}]=x_k[0,b_{N_k}]$, while for $N\in[N_k+1,N_{k+1}]$ we let $\S_{k+1}[a_N,b_N]=\S[a_N,b_N]$.
It will be convenient not to change the enumeration of the orbit segments (except for the first one, which is new), and of the gaps (the first gap of $\S_{k+1}$ coincides with the $N_k$th gap of $\S$). Then $\S_{k+1}$ satisfies $g_N\ge M_{\eps_k}(l_{N+1})$ for all $N\in[N_k,N_{k+1}-1]$ (i.e., for all gaps of $\mathcal S_{k+1}$), hence it can be $\eps_k$-shadowed by the orbit of some point $x_{k+1}\in X$. It is clear that $x_{k+1}$ satisfies~\eqref{szer} with the parameter $k+1$  in place of $k$. This concludes the induction. We let $x_0$ be any accumulation point of the \sq\ $(x_k)_{k\ge1}$. As easily seen, this point satisfies, for all $n\in D$, the inequality 
$$
d(\S(n),T^n(x_0))\le\sum_{j=k_n+1}^\infty\eps_j, 
$$
where $k_n$ is the unique integer $k\ge0$ such that $n\in[a_N,b_N]$ with $N\in[N_k+1,N_{k+1}]$. Since the sums on the right hand side are tails of a convergent series, these distances tend to zero, as claimed.
\end{proof}

\begin{rem}\label{specgen}
It is easily seen that if the domain $D$ of $\S$ in the above lemma has density one then the point $x_0$ from that lemma quasi-generates the same \im s as $\S$.
\end{rem}

The second key lemma requires two rather standard propositions from convex analysis. Although they are well-known to specialists, it is hard to find them in the exact formulation. Thus, we provide them with proofs.  

Let $(\M,d)$ be a compact convex set in a locally convex metric space (the reader may think of $(\M(X),d)$, where $d$ is some standard metric compatible with the weak* topology). The elements of $\M$ will be denoted by the letters $\mu,\nu$.

\begin{prop}\label{ucon} Let $T:\M\to \M$ be a continuous affine transformation. Then the set $\M_T\subset\M$, consisting of $T$-\inv\ elements, is nonempty and for any $\eps>0$ there exists $n_\eps\ge1$ such that, for any $\nu\in\M$ and any $n\ge n_\eps$, we have
$d(\frac1n\sum_{i=0}^{n-1}T^i(\nu),\M_T)<\eps$. 
\end{prop}

\begin{proof}We can assume that $\diam(\M)=1$. Denote $A_n(\nu)=\frac1n\sum_{i=1}^{n-1}T^i(\nu)$. Then,
by convexity of the metric and diameter 1 of $\M$, we easily see that
$$
d(T(A_n(\nu)),A_n(\nu))\le\frac1n.
$$ 
This in turn implies that any limit point of any \sq\ of the form $A_n(\nu_n)$ (with $\nu_n\in\M$) is $T$-\inv. Such limit points exist by compactness, hence we get that $\M_T\neq\emptyset$. Suppose that the second part of the proposition does not hold. This means that there exists $\eps>0$ and an increasing \sq\ $(n_k)_{k\ge1}$ of natural numbers, and a \sq\ $(\nu_k)_{k\ge1}$ of points of $\M$, such that $d(A_{n_k}(\nu_k),\M_T)\ge\eps$ for all $k\ge1$. But we have just proved that all accumulation points of the \sq\ $(A_{n_k}(\nu_k))_{k\ge1}$ belong to $\M_T$, so we have a contradiction. 
\end{proof}

Recall that if $\xi$ is a probability measure on $\M$ then there exists a unique point $\mu\in\M$, called the \emph{barycenter} of $\xi$, such that for every affine continuous function $f$ one has
$$
f(\mu)= \int f(\nu)\,d\mu(\nu).
$$
The barycenter map is denoted by either $\xi\mapsto\mathsf{bar}(\xi)$ or by $\xi\mapsto\int\nu\,d\xi(\nu)$ (the integral in the sense of Pettis).
It is well known that if the set of all probability measures on $\M$ is endowed with the weak* topology then the barycenter map $\xi\mapsto\mathsf{bar}(\xi)$ is continuous.
In the next proposition, the reader may think of $\M$ representing $\M_T(X)$ in a \ds\ \xt, and $\mu$ representing an ergodic measure.
\begin{prop}\label{convlem}
Let $\mu$ be an extreme point of $\M$. Then for any $\eps>0$ there exists $\delta>0$ such that, whenever a probability measure $\xi$ on $\M$
satisfies $d(\mathsf{bar}(\xi),\mu)<\delta$, then 
$$
\xi\{\nu\in\M:d(\mu,\nu)<\eps\}>1-\eps.
$$
\end{prop}
\begin{proof}
If the statement is false then there exists a \sq\ of measures $(\xi_k)_{k\ge1}$ on $\M$ such that $\lim_{k\to\infty}\mathsf{bar}(\xi_k)=\mu$ and $\xi_k\{\nu\in\M:d(\mu,\nu)<\eps\}\le 1-\eps$ for each $k\ge1$. Since the function which associates to a measure $\xi$ the value $\xi(U)$, where $U$ is an open set, is lower semicontinuous in the weak* topology, we get that if $\xi$ is an accumulation point of the \sq\ $(\xi_k)_{k\ge1}$ then
$\xi\{\nu\in\M:d(\mu,\nu)<\eps\}\le 1-\eps$. On the other hand, by continuity of the barycenter map, we have $\mathsf{bar}(\xi)=\mu$.
Since $\mu$ is an extreme point of $\M$, the only measure on $\M$ with barycenter at $\mu$ is the Dirac measure $\delta_{\mu}$. We conclude that $\xi=\delta_{\mu}$. This is a contradiction, since $\delta_{\mu}\{\nu\in\M:d(\mu,\nu)<\eps\}=1$.
\end{proof}


We proceed with the second key lemma needed in the proof of Theorem \ref{lgp}.

\begin{lem}\label{qg}
Let $\mu$ be an ergodic measure on a \tl\ \ds\ \xt\ which has the weak specification property. Let $x_0\in X$ be quasi-generic for $\mu$ and let $\mathcal{J}=(n_k)_{k\ge1}$ be a \sq\ along which $x_0$ generates $\mu$. Then there exists a point $\bar x_0\in X$ generic for $\mu$ and a set $\mathbb M\subset\N$ of upper density one achieved along a sub\sq\ of $\mathcal{J}$, such that 
$$
\lim_{n\in\mathbb M}\ d(T^n\bar x_0,T^nx_0)=0.
$$
\end{lem}

\begin{proof} We start by fixing a summable \sq\ of positive numbers $(\eps_k)_{k\ge1}$. 
In view of Lemma \ref{goodshad} and Remark \ref{specgen},
it suffices to construct a specification $\S$ on a domain $D$ satisfying the following four conditions: 
\begin{enumerate}
	\item the assumptions of Lemma \ref{goodshad}, 
	\item the domain $D$ of $\S$ has density one,
	\item $\lim_{n\in\mathbb M} d(\S(n),T^n(x_0))=0$, where $\mathbb M\subset D$ has upper density one achieved along a sub\sq\ of $\mathcal{J}$,
	\item $\S$ is generic for $\mu$.
\end{enumerate}

We choose a \sq\ of positive integers $(l_k)_{k\ge0}$. The \sq\ should grow so fast that the ratios $\frac{M_{\eps_k}(l_k)}{l_k}$ are all smaller than $1$ and tend to zero. For each $k\ge1$ we let $L_k=l_k+M_{\eps_k(l_k)}$. 
Next, we replace $\mathcal J$ by a fast growing sub\sq\ and from now on $\mathcal J=(n_k)_{k\ge1}$ will denote that sub\sq. Initially we require that the ratios $\frac{l_k}{n_k}$ and $\frac{n_k}{n_{k+1}}$ tend to zero as $k$ grows. More  conditions on the speed of growth of the \sq s $(l_k)_{k\ge1}$ and $(n_k)_{k\ge1}$ will be specified later. 

The specification $\S$ is created in three steps. The first auxiliary specification $\S'$ is just a partition of the orbit of $x_0$ without gaps. We begin by partitioning it into segments of length $L_1$ until we cover the coordinate $n_1$. Then we continue by partitioning the remaining part of the orbit of $x_0$ into segments of length $L_2$ until we cover the coordinate $n_2$ and so on. To be precise, we create segments $\S'[a_N,b_N]=x_0[a_N,b_N]$ (where $N\ge1$) satisfying: 
\begin{enumerate}[(i)]
	\item $a_1=0$, 
	\item for $N\ge2$, $a_N=b_{N-1}+1$, 
	\item $b_N=a_N+L_k-1$, for $N\in[N_{k-1}+1,N_k]$, where
	\item for each $k\ge1$, $N_k$ is such that $n_k\in[a_{N_k},b_{N_k}]$
	\end{enumerate}
(for consistency of the notation, we have let $N_0=0$). It is elementary to see that, since the ratios $\frac{M_{\eps_k}(l_k)}{l_k}$ and $\frac{l_k}{n_k}$ tend to zero, the ratios $\frac{b_{N_k}}{n_k}$ tend to $1$. Thus, by Remark \ref{eqseq}, the point $x_0$ generates $\mu$ along the \sq\ $(b_{N_k})_{k\ge1}$. From now on, we redefine the \sq\ $\mathcal J$ to be $(b_{N_k})_{k\ge1}$ (and let $n_k=b_{N_k}$; we also let $n_0=0$). This new \sq\ still satisfies $\frac{n_k}{n_{k+1}}\to0$. 

The empirical measures $\mu_{x_0[0,n_k]}$ tend to $\mu$. Since $n_{k-1}$ is eventually negligible in comparison with $n_k$, the following holds:
\begin{equation}\label{convv}
\text{the empirical measures $\mu_{x_0[n_{k-1}+1,n_k]}$ tend to $\mu$ as $k$ grows}.
\end{equation}

The second auxiliary specification $\S''$ is obtained from $\S'$ by truncating all orbit segments (except the first one) on the left, to allow for future shadowing. More precisely, we let $a'_1=a_1=0$ and for any $k\ge1$ and any $N\in[N_{k-1}+1,N_k]$ (except for $N=1$) we let $a'_N=a_N+M_{\eps_k}(l_k)$ (since $M_{\eps_k}(l_k)<l_k$, we have $a_N'<b_N$). Then, on the new domain
$$
D=\bigcup_{N\ge1}[a'_N,b_N],
$$
we define the specification $\S''$ by $S''[a'_N,b_N]=x_0[a'_N,b_N]$. This new specification has, for $N\in[N_{k-1}+1,N_k]$ (except for $N=1$), orbit segments of length $l_k$ preceded by gaps of size $M_{\eps_k}(l_k)$ (the first orbit segment has length $L_1$ and no preceding gap).


It should be quite obvious that the lower density of $D$ is achieved along the \sq\ $b_{N_k}+M_{\eps_{k+1}(l_{k+1})}$
(this is the end of the first gap larger than all preceding gaps). Because the ratios $\frac{M_{\eps_k}(l_k)}{l_k}$ tend to zero, by choosing the numbers $n_k$ (and hence $b_{N_k}$) sufficiently large in comparison with $M_{\eps_{k+1}}(l_{k+1})$, we can arrange that the density of $D$ equals one, as required in (2).

Now, we have to go back to the choice of the \sq s $(l_k)$ and $(n_k)$ and impose more conditions on the speed of their growth. We select numbers $\delta_k\le\eps_k$ according to Proposition \ref{convlem} with respect to the numbers $\eps_k$ and the ergodic measure $\mu$ in the role of the extreme point of the compact convex set $\M_T(X)$. If the numbers $l_k$ are (a priori) chosen large enough, using Proposition \ref{ucon}, we can arrange that 
\begin{itemize}
\item the empirical measures $\mu_{x_0[a'_N,b_N]}$ with $N\in[N_{k-1}+1,N_k]$ are $\frac{\delta_k}3$-close to some \im s henceforth denoted by $\mu_N$. 
\end{itemize}
Also, by imposing fast enough growth of the numbers $n_k$, we may achieve that 
\begin{itemize}
	\item the empirical measure $\mu_{x_0[n_{k-1}+1,n_k]}$ is $\frac{\delta_k}3$-close to $\mu$ (see \eqref{convv}),
	\item the empirical measures $\mu_{x_0[a_N,b_N]}$ with $N\in[N_{k-1}+1,N_k]$ are $\frac{\delta_k}3$-close to the respective empirical measures $\mu_{x_0[a'_N,b_N]}$ (and hence $\frac23\delta_k$-close to $\mu_N$).
\end{itemize}

Clearly,  the empirical measure $\mu_{x_0[n_{k-1}+1,n_k]}=\mu_{x_0[a_{N_{k-1}+1},b_{N_k}]}$ equals the arithmetic average of the measures $\mu_{x_0[a_N,b_N]}$ with $N\in[N_{k-1}+1,N_k]$. By convexity of the metric, $\mu$ is $\delta_k$-close to the arithmetic average of the \im s $\mu_N$ with $N\in[N_{k-1}+1,N_k]$. By Proposition \ref{convlem}, vast majority of the \im s $\mu_N$ are $\eps_k$-close to $\mu$, and hence the corresponding empirical measures $\mu_{x_0[a'_N,b_N]}$ are $2\eps_k$-close to $\mu$ (we are using $\delta_k<\eps_k$). More precisely, there are less than $\eps_k(N_k-N_{k-1})$ parameters $N\in[N_{k-1}+1,N_k]$ (we will call them ``bad''), for which the measure $\mu_{x_0[a'_N,b_N]}$ is not $2\eps_k$-close to $\mu$. 

We can now perform the third step in creating the specification $\S$. This is done by replacing in $\S''$ the segments $x_0[a'_N,b_N]$ corresponding to ``bad'' parameters $N\in[N_{k-1},N_k]$ by orbit segments (of the same length) whose associated measures are $2\eps_k$-close to $\mu$. For example, we can choose one ``good'' parameter $N$ (there exists such an $N$) and use the corresponding segment $x_0[a'_N,b_N]$ everywhere we need to make a replacement. This concludes the construction of~$\S$.

We need to verify that $\S$ satisfies the desired 4 properties. And so:
\begin{enumerate}
	\item It is clear that the specification $\S$ was created in accordance with the assumptions of Lemma \ref{goodshad}.
	\item As we have already remarked, the domain $D$ has density one.
	\item Note that $\S''$ agrees with the orbit of $x_0$ on $D$ (which has density one). Then $\S$ differs from $\S''$ on a set whose frequency in the interval $[n_{k-1}+1,n_k]$ is at most $\eps_k$. Thus the set of the integers $n$ for which $\S(n)\neq\S''(n)$ (or $\S(n)$ is not defined) has lower density zero achieved along the \sq\ $\mathcal J$.
The complementary set $\mathbb M$ has upper density one achieved along $\mathcal J$ and on this set we have $\S(n)=T^n(x_0)$ (which trivially implies the required condition $\lim_{n\in\mathbb M}\ d(\S(n),T^nx_0)=0$).
	\item Consider a long initial segment of $\S$, say $\S|_{[1,n]\cap D}$, and let $k$ be such that $N\in[N_{k-1}+1,N_k]$, where $N$ is determined by the inclusion $n\in[a_N,b_N]$. Then $\S|_{[1,n]\cap D}$ consists essentially of segments of two lengths: $l_{k-1}$, whose associated empirical measures are $2\eps_{k-1}$-close to $\mu$, and $l_k$, whose associated empirical measures are $2\eps_k$-close to $\mu$ (in either case we have $2\eps_{k-1}$-closeness). This closeness need not apply to the initial part left of the coordinate $n_{k-1}$, and to the terminal, perhaps incomplete, orbit segment whose length does not exceed $L_k$. Since both $n_{k-1}$ and $L_k$ is negligible in comparison with $n_k$ (and hence with $n$), the two extreme pieces can be ignored and we get that the empirical measure associated with $\S|_{[1,n]\cap D}$ is (nearly) $2\eps_{k-1}$-close to $\mu$. Since $k$ tends to infinity as $n$ grows, $\S$ is generic for $\mu$.
\end{enumerate}
\end{proof}

\begin{lem}\label{new}
Let \xt\ be a \tl\ \ds\ and let $\mu$ be an \im\ on $X$. 
For each $\eps>0$ there exists $\delta>0$ which satisfies the following:

Let $\P$ be a partition of $X$ whose all atoms have diameter not exceeding $\delta$. 
Let $\S$ be a finite specification consisting of $N_1$ orbit segments of lengths $l$ 
separated by some gaps ($N_1$ and $l$ are arbitrary natural numbers, the gaps are also arbitrary):
$$
\S[a_N,a_N+l-1]=x_N[0,l-1],
$$ 
where $a_1\ge0$ and, for each $N\in[1,N_1]$, we have $x_N\in X$ and $a_{N+1}-a_N\ge l$. Assume that for each $B\in\P^l=\bigvee_{n=0}^{l-1}T^{-n}(\P)$ the frequency \emph{relative to $(a_N)_{N\in[1,N_1]}$}:
$$
\frac{|\{N\in[1,N_1]: x_N\in B\}|}{N_1}
$$
is $\frac\delta{|\P^l|}$-close to $\mu(B)$ (this imposes that $N_1$ must in fact be huge). Then the empirical measure associated with $\S$,
$$
\mu_\S=\frac1{|D|}\sum_{n\in D}\delta_{\S(n)},
$$ 
where $D=\bigcup_{N=1}^{N_1}[a_N,a_N+l-1]$ is the domain of $\S$, is $\eps$-close to $\mu$. 
\end{lem}

\begin{proof} Regardless of what metric $d$ compatible with the weak* topology on $\M(X)$ we are using, there exists a finite family of continuous $[0,1]$-valued functions, say $f_1,\dots,f_K$, and a small positive number $\gamma$ such that if 
$$
\left|\int f_k\,d\mu_1 - \int f_k\,d\mu_2\right|<3\gamma   
$$
for each $k\in[1,K]$, then $d(\mu_1,\mu_2)<\eps$. Further, there exists $\beta$ such that each of the finitely functions $f_k$ varies on each $\beta$-ball in $X$ by less than $\gamma$. We let $\delta=\min\{\beta,\gamma\}$. Let $\P$ be a partition of $X$ as in the formulation of the lemma. Observe that if we replace each of the functions $f_k$ by a function $\bar f_k$ constant on the atoms of $\P$ (say, assuming on each atom the supremum of $f_k$ over that atom), then the integral of $\bar f_k$ with respect to any probability measure differs from the integral of $f_k$ by at most~$\gamma$. So, in order to show that $d(\mu_\S,\mu)<\eps$, it suffices to show that
$$
\left|\int f\,d\mu_\S - \int f\,d\mu\right|<\gamma,   
$$
for any $[0,1]$-valued (not necessarily continuous) function $f$ constant on the atoms of $\P$. 
For such a function $f$ we have
$$
\int f\,d\mu_\S = \frac1{|D|}\sum_{N=1}^{N_1} \sum_{n=0}^{l-1}f(T^nx_N)=
\frac1{N_1}\sum_{N=1}^{N_1}\frac1l\sum_{n=0}^{l-1}f(T^nx_N)
$$
(we are using the obvious fact that $|D|=N_1l$). Note that if, for some $N,N'\in[1,N_1]$, teh points $x_N$ and $x_{N'}$ belong to the same atom $B$ of $\P^l$ then the averages $\frac1l\sum_{n=0}^{l-1}f(T^nx_N)$ and $\frac1l\sum_{n=0}^{l-1}f(T^nx_{N'})$ are equal, so, we can replace them by $\frac1l\sum_{n=0}^{l-1}f(T^nx_B)$, where $x_B$ is a point in $B$ not depending on $N$. Then, our integral becomes
$$
\int f\,d\mu_\S = \sum_{B\in\P^l}\frac{|\{N\in[1,N_1]:x_N\in B\}|}{N_1}\frac1l\sum_{n=0}^{l-1}f(T^nx_B).
$$
By assumption, the coefficient $\frac{|\{N\in[1,N_1]:x_N\in B\}|}{N_1}$ equals $\mu(B)$ up to $\frac\delta{|\P^l|}$, all the more up to $\frac\gamma{|\P^l|}$. Since the averages $\frac1l\sum_{n=0}^{l-1}f(T^nx_B)$ do not exceed $1$ we obtain that $\int f\,d\mu_\S$ equals
$$
\sum_{B\in\P^l}\mu(B)\frac1l\sum_{n=0}^{l-1}f(T^nx_B)
$$
up to $\gamma$.
Finally, observe that the latter sum equals $\int \frac1l\sum_{n=0}^{l-1}f\circ T^n\,d\mu$, which, by invariance of $\mu$, equals $\int f\,d\mu$. We have shown that $|\int f\,d\mu_\S - \int f\,d\mu|<\gamma$, as needed.
\end{proof}

The next proposition is our last preparatory fact before the proof of Theorem~\ref{lgp} It is also a standard fact (this time from ergodic theory), whose exact formulation is hard to find. Thus, we provide it with a proof.

\begin{prop}\label{missing}
Let \xt\ be a \tl\ \ds. Let $x$ be a point quasi-generic for an ergodic measure $\mu$ and let $\mathcal J=(n_k)_{k\ge1}$ be a \sq\ along which $x$ generates $\mu$. Fix a positive integer $L$. Then there exist two increasing \sq s\ of positive integers: $(a_N)_{N\ge0}$ and $(N_k)_{k\ge1}$ satisfying the following conditions:
\begin{enumerate}
	\item for each $N\ge1$ the difference $a_{N+1}-a_N$ equals either 
	$L$ or $L+1$, 
	\item $\lim_{k\to\infty}\frac{a_{N_k}}{n_k}=1$ (i.e., the \sq s $(n_k)_{k\ge1}$ and $(a_{N_k})_{k\ge1}$ are equivalent),
	\item $x$ generates $\mu$ \emph{relatively w.r.t.\ the \sq\ $(a_N)_{n\ge1}$}, along $(N_k)_{k\ge1}$, i.e., 
	$$
	\lim_{k\ge1}\frac1{N_k}\sum_{N=1}^{N_k}\delta_{T^{a_N}(x)}=d\mu.
	$$
\end{enumerate} 
\end{prop}

\begin{proof}
There exists an ergodic measure-preserving system $(Y,\nu,S)$ disjoint from $(X,\mu,T)$ (in the sense of Furstenberg).\footnote{Two measure preserving systems are \emph{disjoint} if their only joining is their product. An example of a system disjoint from $(X,\mu,T)$ is an irrational rotation by $e^{2\pi it}$, where $t$ is rationally independent from all numbers $s$ such that $e^{2\pi is}$ is an eigenvalue of $(X,\mu,T)$ (there are at most countably many values to be avoided).} By a standard application of Rokhlin towers, there exists a set $A$ visited by $\nu$-almost every orbit in $Y$ infinitely many times with only two gap sizes between consecutive visits, $L$ and $L+1$. There exists a \tl\ model of $(Y,\nu,S)$ in which the set $A$ is clopen, i.e., its indicator function, denoted by $F$, is continuous. Let $y\in Y$ be a point generic for $\nu$ and let $(a_N)_{N\ge1}$ denote the \sq\ of times of visits of the orbit of $y$ in $A$ (this \sq\ has only two gap sizes: $L$ and $L+1$, as required in (1)). The pair $(x,y)$ generates, along the \sq\ $\mathcal J$, some joining of $\mu$ and $\nu$. By disjointness, this joining equals the product measure $\mu\times\nu$ on $X\times Y$. This implies that,
for any continuous function $f$ on $X$, we have
\begin{align*}
\lim_{k\to\infty}\frac1{n_k}\sum_{n=1}^{n_k}f(T^nx)F(S^ny)&=\int f\,d\mu\cdot\nu(A), \text{ \ and}\\
\lim_{k\to\infty}\frac1{n_k}\sum_{n=1}^{n_k}F(S^ny)&=\nu(A).
\end{align*}
Given $k\ge1$, let $N_k$ denote the largest $N$ such that $a_N\le n_k$. Observe that since $(a_N)_{N\ge1}$ has bounded gaps, while $(n_k)_{k\ge1}$ tends to infinity, the ratios $\frac{a_{N_k}}{n_k}$ tend to~$1$, as required in (2).
Since $F(S^ny)=1$ if and only if $n=a_N$ for some $N$ (otherwise $F(S^ny)=0$), we can rewrite the above limits as
\begin{align*}
\lim_{k\to\infty}\frac1{n_k}\sum_{N=1}^{N_k}f(T^{a_N}x)&=\int f\,d\mu\cdot\nu(A), \text{ \ and}\\
\lim_{k\to\infty}\frac{N_k}{n_k}&=\nu(A).
\end{align*}
Dividing sidewise, we get 
$$
\lim_{k\to\infty}\frac1{N_k}\sum_{N=1}^{N_k}f(T^{a_N}x)=\int f\,d\mu.
$$
Since this is true for any continuous function $f$ on $X$, we have proved (3).
\end{proof}

\section{The main proof}\label{S4}
\begin{proof}[Proof of Theorem \ref{lgp}] Let $\mathcal J=(n_k)_{k\ge1}$ be a \sq\ along which $y$ generates $\nu$.
It suffices to construct a point $x_0$ such that the pair $(x_0,y)$ generates $\xi$ along a sub\sq\ $\mathcal J'$ of $\mathcal{J}$. Clearly, such an $x_0$ generates $\mu$ along $\mathcal J'$ and, by Lemma~\ref{qg}, there will then exist a point $x$ generic for $\mu$ and such that 
$$
\lim_{n\in\mathbb M}\ d(T^nx,T^nx_0)=0,
$$
where $\mathbb M$ is a set of upper density one achieved along a sub\sq\ $\mathcal J''$ of $\mathcal J'$. Note that then 
the pair $(x,y)$ still generates $\xi$ along $\mathcal J''$, so the proof will be completed. 
\smallskip

We fix a summable \sq\ of positive numbers $(\eps_k)_{k\ge1}$ and an increasing \sq\ of natural numbers $l_k$, such that 
\begin{equation}\label{eq}
\lim_{k\to\infty}\frac{M_{\eps_k}(l_k)}{l_k}=0.
\end{equation}

Next, we let $(\P_k)_{k\ge1}$ be a \sq\ of measurable partitions of $X$ such that, for each $k\ge1$, the diameters of the atoms of $\P_k$ do not exceed the number $\delta_k$ obtained from Lemma \ref{new} for the measure $\mu$ and $\eps_k$ in the role of $\eps$.  

The atoms of the partitions $\P_k^l=\bigvee_{i=0}^{l-1}T^{-i}(\P_k)$, where $k\ge1$ and $l\ge1$, will be called \emph{blocks} of $X$, while the atoms $\P_k^{l_k}$ will be called \emph{blocks of order $k$ of $X$}. 

Likewise, we let $(\Q_k)_{k\ge1}$ be a \sq\ of partitions of $Y$ with diameters bounded by $\delta_k$. 
We can easily arrange the partitions $\Q_k$ so that the orbit of $y$ avoids the boundaries of the atoms of $\Q_k$ for each $k\ge1$. The atoms of $\Q_k^l=\bigvee_{i=0}^{l-1}S^{-i}(\Q_k)$, where $k\ge1$ and $l\ge1$, will be called \emph{blocks} of $Y$ and the atoms of $\Q_k^{l_k}$ will be called \emph{atoms of order $k$ of $Y$}. 

Note that if we apply the maximum metric in $X\times Y$ then the rectangular atoms of the partitions $\P_k\otimes\Q_k$ have diameters bounded by $\delta_k$ as well. We now choose some very small positive numbers $\gamma_k$ so that, for each $k\ge1$ we have 
$$
2\gamma_k+\gamma_k^2<\frac{\delta_k}{|\P_k^l\otimes\Q_k^l|}. 
$$

Because $y$ is generic for $\nu$ along $\mathcal J$, and its orbit avoids the boundaries of the blocks, the orbit of $y$ visits each block $C$ of $Y$ with the frequency evaluated at times $n_k$ converging to $\nu(C)$. 

Using successively Proposition \ref{missing} with the parameters $L_k=l_k+M_{\eps_k}(l_k)$ in the role of $L$, and replacing, if necessary, the \sq\ $\mathcal J$ by a rapidly growing sub\sq\ $\mathcal J'$ (from now on $(n_k)_{k\ge1}$ will denote $\mathcal J'$), we can arrange two increasing \sq s of positive integers: $(a_N)_{N\ge1}$ and $(N_k)_{k\ge0}$ starting with $N_0=0$, satisfying the following conditions:
\begin{enumerate}
	\item $\displaystyle{\lim_{k\to\infty}\frac{a_{N_k}}{n_k}=1}$,
	\item for each $k\ge1$ and each $N\in[N_{k-1},N_k-1]$ the difference $a_{N+1}-a_N$ equals either $L_k$ or $L_k+1$,\footnote{The lemma, as it is stated, does not allow to ensure that the gap $a_{N_k}-a_{N_k-1}$ (the first gap following the series of gaps of sizes $L_k$ or $L_k+1$) equals either $L_{k+1}$ or $L_{k+1}+1$. A priori, it may come out too small. However, replacing the set $A$ in the proof of that lemma by $T^{-j}A$ and enlarging, if necessary, $n_k$, we can shift the term $a_{N_k}$ by an arbitrary positive integer $j\le l_{k+1}$ to the right, and in this manner adjust the gap.}
	\item if we denote by $C_N$ the unique block of order $k$ of $Y$ containing $S^{a_N}y$, then, for any 
	block $C$ of order $k$ of $Y$, one has
	$$
	\left|\frac1{N_k-N_{k-1}}|\{N\in[N_{k-1},N_k-1]:C_N=C\}|-\nu(C)\right|<\gamma_k.
	$$
\item if, in addition, $C$ satisfies $\nu(C)>0$, then also 
$$
|\{N\in[N_{k-1},N_k-1]:C_N=C\}|>\frac1{\gamma_k}.
$$
\end{enumerate}
The condition~(1) says that $\mathcal J'$ and $\tilde{\mathcal J}'=(a_{N_k})_{k\ge1}$ are equivalent. In particular $y$ generates $\nu$ along the \sq\ $(a_{N_k})_{k\ge1}$ and if we find $x_0$ using $\tilde{\mathcal J}'$, the same $x_0$ will serve for $\mathcal J'$.

(*) Fix some $k\ge1$. Let $\xi(B|C)=\frac{\xi(B\times C)}{\nu(C)}$, where $B$ and $C$ are blocks of order $k$ of $X$ and $Y$, respectively, with $\nu(C)>0$. 
For every such $C$, the numbers $\xi(B|C)$, with $B$ ranging over all blocks of order $k$ of $X$, form a probability vector. By~(4), this vector can be approximated up to $\gamma_k$ (at each coordinate) by a rational probability vector with entries 
$$
\frac{r(B,C)}{|\{N\in[N_{k-1},N_k-1]:C_N=C\}|},
$$ 
where each $r(B,C)$ is a nonnegative integer. We can thus create a finite \sq\ $(B_N)_{N\in[N_{k-1},N_k-1]}$ of blocks of order $k$ of $X$, so that, for every pair of blocks $B,C$ of order $k$ in $X$ and $Y$, respectively, we have
$$
|\{N\in[N_{k-1},N_k-1]:C_N=C\text{ and }B_N=B\}|=r(B,C).
$$
Then, for each pair $B,C$ as above, with $\nu(C)>0$, we have
\begin{multline*}
\frac{r(B,C)}{N_k-N_{k-1}}=\\
\frac{r(B,C)}{|\{N\in[N_{k-1},N_k-1]:C_N=C\}|}\cdot \frac{|\{N\in[N_{k-1},N_k-1]:C_N=C\}|}{N_k-N_{k-1}},
\end{multline*}
where (by the choice of the integers $r(B,C)$) the first fraction equals $\xi(B|C)$ up to~$\gamma_k$, and, by (3), the second fraction equals $\nu(C)$, also up to $\gamma_k$. So, $\frac{r(B,C)}{N_k-N_{k-1}}$ equals $\xi(B\times C)$ up to $2\gamma_k+\gamma_k^2$, which is less than $\frac{\delta_k}{|\P_k^l\otimes\Q_k^l|}$.

Now, we create a finite specification $\bar{\S}_k$ in $X\times Y$, as follows. For each $N\in[N_{k-1},N_k-1]$ we choose a point $x_N\in B_N$ and we let 
$$
\bar{\S}_k[a_N,a_N+L_k-1]=(x_N,S^{a_N}y)[0,L_k-1]
$$
(the starting point of the $N$th orbit segment falls in $(B_N,C_N)$, the second coordinate agrees, along the entire specification, with the orbit of $y$).
Note that by (2), the gaps in the domain of $\bar{\S}_k$ have only two sizes, $0$ or $1$. Lemma \ref{new} now guarantees that the empirical measure $\mu_{\bar{\S}_k}$ is $\eps_k$-close to $\xi$. 

Let $\bar\S$ be the infinite specification in $X\times Y$ defined as follows: for each $k\ge1$ and each $N\in[N_{k-1},N_k-1]$ we let
$$
\bar{\S}[a_N+M_{\eps_k}(l_k),a_N+L_k-1]=\bar{\S}_k[a_N+M_{\eps_k}(l_k),a_N+L_k-1].
$$
It is fairly obvious that $\bar{\S}$ generates $\xi$ along the \sq\ $\tilde{\mathcal J}'$ (by \eqref{eq}, the fact that of the intervals of the domain are slightly trimmed on the left does not affect the convergence).

Let us denote by $\S$ the projection of $\bar\S$ to the first coordinate. This infinite specification in $X$ satisfies all requirements of Lemma \ref{goodshad}. That lemma allows to find a point $x_0$ whose orbit shadows the specification $\S$ with an increasing accuracy. Clearly, the pair $(x_0,y)$ shadows $\bar\S$ equally well. It is also clear that the domain of $\bar\S$ has density $1$, which (by Remark \ref{specgen}) implies that the pair $(x_0,y)$ generates $\xi$ along $\tilde{\mathcal J}'$, and hence also along $\mathcal J'$. We have achieved all that was necessary to complete the proof.
\end{proof}

\begin{rem}
It is possible to modify the proof and avoid the use of Lemma~\ref{qg} (see below). Although the main proof itself becomes slightly longer, one can skip that lemma and the two auxiliary propositions altogether. There are two reasons why we have decided to present the longer argument: 
\begin{enumerate}
	\item Lemma \ref{qg} is a generalization of Kamae's Theorem~1 in~\cite{K} and has some independent value of its own. It may turn out useful in further studies of systems with weak specification. 
	\item By following the framework of the original proof, we show that T.\ Kamae has insightfully laid ground for further generalizations.
\end{enumerate}
\end{rem}

\begin{proof}[Sketch of the modified proof]
Go to the paragraph marked by (*). Divide the blocks $B$ (of order $k$ of $X$) into two families: $\B$, of those whose associated empirical measures are close to $\mu$ and the rest.
By the mean ergodic theorem, for large enough $k$, the joint measure of the blocks in $\B$ is very close to $1$. Thus, by an insignificant renormalization, we can make the vector of conditional probabilites $\xi(B|C)$, with $B$ ranging over $\B$ (and $C$ fixed) probabilistic. From here we proceed as it is described except that each time we refer to $B$ we use only the blocks from~$\B$. The specification $\bar\S_k$ will then have its $X$-coordinate consisting exclusively of blocks $B\in\B$. The specification $\S$ (the $X$-projection of $\bar\S$) will consist of blocks whose empirical measures are getting closer and closer to $\mu$. If the numbers $a_{N_k}$ grow sufficiently fast in comparison to the lengths $l_k$ then, by an identical argument as in the proof of Lemma \ref{qg}, the specification $\S$ will be generic for $\mu$ and so will be the point $x_0$ shadowing $\S$. Lemma \ref{qg} becomes irrelevant.
\end{proof}

\end{document}